\documentclass[12pt,oneside,english]{amsart}
\usepackage[latin1]{inputenc}
\usepackage{amssymb}

\makeatletter
\newtheorem{theorem}{Theorem}

\newtheorem{example}{Example}

\usepackage{babel}

\makeatother
\begin{document}
\baselineskip=17pt

\title[On the Newman sum ]{On the Newman sum over multiples of a prime
with a primitive or semiprimitive root~2}

\author{Vladimir Shevelev}
\address{Departments of Mathematics \\Ben-Gurion University of the
 Negev\\Beer-Sheva 84105, Israel. e-mail:shevelev@bgu.ac.il}


\begin{abstract}
We prove that if  $S_p(x)$ is the Newman sum over p-multiples where
$p$  is a prime with a primitive (semiprimitive) root 2 then
$S_p(2^p x)=pS_p(2x)\;\;(S_p(2^p x)=(-1)^{\frac{p-1}{2}}pS_p(2x))$.
We consider the case of $p=17$ as well.
\end{abstract}

\maketitle
\section{Introduction}
Consider  for $x,m, l\in\mathbb{N},\;\;l\in[0,m-1]$,  the Newman sum

$$
S_m(x)=S_{m,0}(x)=\sum_{0\leq n <x,\;n\equiv 0(mod
m)}(-1)^{\sigma(n)},
$$

where $\sigma(n)$ is the number of 1's in the binary expansion of
$n$.

In \cite{6} we proved that if 2 is a primitive root of a prime $p$
then $S_p(2^p)=p$. Now we prove a more general relation.

\begin{theorem}\label{t1}
If 2 is a primitive root of a prime $p$ then for any natural $x$
\begin{equation}\label{1}
S_p(2^px)=pS_p(2x).
\end{equation}
\end{theorem}

Furthermore, let 2 be not a primitive root of an odd prime $p$. We
say that 2 is a  \slshape semiprimitive  root \upshape of $p$ if 2
has the order $\frac{p-1}{2}$ modulo $p$ and the congruence
$2^x\equiv -1\pmod{p}$ is not solvable.

\begin{example}
2 has the order $\frac{p-1}{2}$ for $p=7,17,23$ but only for $p=17$
the congruence $2^x\equiv -1(mod 17)$ has a solution $(x=4)$.
Therefore, by the definition, 2 is a semiprimitive root for $7,23$
(but not  for 17).
\end{example}
    Note that, if 2 is a semiprimitive root of $p$ then for every
integer $a\in[1,p-1]$ there exist $j=j_a$ such that either $a\equiv
2^j\pmod{p}$ or $a\equiv -2^j\pmod{p}$.

Our statement similar to Theorem 1 in the case of 2 is a
semiprimitive root of $p$ is the following.
\newpage

\begin{theorem}\label{t2}

If 2 is a semiprimitive root of a prime $p$ then for any natural $x$

\begin{equation}\label{2}
S_p(2^px)=(-1)^{\frac{p-1}{2}}pS_p(2x).
\end{equation}
\end{theorem}

Theorems \ref{t1} and \ref{t2} show that if 2 is primitive or
semiprimitive root of an odd prime $p$ then

\begin{equation}\label{3}
S_p(x)=O\left(x^\frac{\ln p}{(p-1)\ln 2}\right).
\end{equation}

and open a way similar to \cite{4} to get the sharp estimates for
$S_p(x)$ in considered cases , i.e. to generalize  the Coquet's
theorem (see\cite{1}, p.98-99).

On the other hand, (\ref{3}) makes more precise the remainder term
of the Gelfond theorem in the considered case:

$$
G_p^{(i)}(x)=\sum_{0\leq n < x,n\equiv 0(mod p),\sigma(n)\equiv
i(mod 2)}1=\frac{x}{2p}+O\left(x^{\frac{\ln p}{(p-1)\ln
2}}\right),i=0,1,
$$

instead of $O\left(x^{\frac{\ln 3}{\ln 4}}\right)$ in
\cite{3}.Moreover, these estimates are unimprovable.

\section{Proof of Theorem 1}
We again use the formula (cf.\cite{2})

\begin{equation}\label{4}
S_p(2^k)= \frac 1 p
\sum^{p-1}_{l=1}\prod^{k-1}_{j=0}(1-\omega_p^{l2^j}),
\end{equation}

where $\omega_p\neq 1$ is a primitive root of $1$ of the power $p$.
By (\ref{4}) we have also for $k\geq p$

\begin{equation}\label{5}
S_p(2^{k-p+1})= \frac 1 p
\sum^{p-1}_{l=1}\prod^{k-p}_{j=0}(1-\omega_p^{l2^j}).
\end{equation}

Let us consider the quotient
$$
Q=\frac{\prod^{k-1}_{j=0}\left(1-\omega_p^{l2^j}\right)}{\prod^{k-p}_{j=0}
\left(1-\omega_p^{l2^j}\right)}=
\prod^{k-1}_{k-(p-1)}\left(1-\omega_p^{l2^j}\right),\;\;l=1,2,\ldots,p-1.
$$

Considering here the substitutions $j-k+p=t,\;l_1=l\cdot 2^{k-p}$ we
find
\newpage

\begin{equation}\label{6}
Q=\prod^{p-1}_{t=1}\left(1-\omega_p^{l,2^t}\right)
\end{equation}

Since 2 is a primitive root of $p$ then independently on $l$

\begin{equation}\label{7}
Q=\prod^{p-1}_{t=1}\left(1-\omega_p^t\right)
\end{equation}

Note that

$$
\prod^{p-1}_{t=1}\left(x-\omega_p^t\right)=\frac{x^p-1}{x-1}=1+x+\ldots+x^{p-1}.
$$

Therefore, by (\ref{7}) $Q=p$ and according to (\ref{4})-(\ref{5})

$$
S_p\left(2^p\cdot 2^{k-p}\right)=pS_p\left(2\cdot 2^{k-p}\right).
$$

Thus, for $x=2^n,\;n\geq 0$, we have

\begin{equation}\label{8}
S_p(2^px)=p S_p(2x).
\end{equation}

Finally, using the additive properties of $S_p$ as in \cite{4} we
obtain (\ref{8}) for any nonnegative integer $x$. $\blacksquare$

As well, note that if $S_p([x,y))$ denotes the difference
$S_p(y)-S_p(x)$ then we have

\begin{equation}\label{9}
S_p([2^p x, 2^p y))=p S_p([2x, 2y)).
\end{equation}

\section{Proof of Theorem 2.}

Since 2 is a semiprimitive root of $p$ then instead of (\ref{6})
independently of $l_1$ we have

$$
Q=\prod^{\frac{p-1}{2}}_{j=1}\left(1-\omega_p^{2^j}\right)^2=(-1)^{\frac{p-1}{2}}
\prod^{\frac{p-1}{2}}_{j=1}\left(1-\omega_p^{2^j}\right)
\prod^{\frac{p-1}{2}}_{j=1}\left(\omega_p^{2^j}-1\right)=
$$
$$
=(-1)^{\frac{p-1}{2}}\omega^{2+2^2+\ldots+2^{\frac{p-1}{2}}}
\prod^{\frac{p-1}{2}}_{j=1}\left(\left(1-\omega_p^{2^j}\right)\left(1-\omega_p^{-2^j}\right)\right)=
$$
$$
=(-1)^{\frac{p-1}{2}}\prod^{p-1}_{k=1}\left(1-\omega_p^k\right)=(-1)^{\frac{p-1}{2}}p.
$$
\newpage

Now in this case similar to (\ref{8}) we obtain (\ref{2}) and the
following relation

$$
S_p([2^p x, 2^p y))=(-1)^{\frac{p-1}{2}}p S_p([2x, 2y)).\blacksquare
$$

\section{Case of $p=17$}

Now we give a relation for the first number of the Drmota-Skalba
primes \cite{2} $p=17$ for which 2 is neither primitive nor
semiprimitive root.

\begin{theorem}\label{t3}
$$
S_{17}(2^{17} x)= 34 S_{17}(2^9 x)-17 S_{17}(2x),\;\;x\in\mathbb{N}.
$$
\end{theorem}

In particular, in the case of $x=1$ we have
$$
S_{17}(2^{17})= 34 S_{17}(2^9)-17 S_{17}(2)=31\cdot 21-17=697.
$$

Using Theorem \ref{3} as in \cite{4} it could be proved that

$$
S_{17}(x)=O(x^\alpha)
$$

with $\alpha=\frac{\ln(17+4\sqrt{17})}{\ln 256}=0.633220353\ldots$.
It is essentially more than $\frac{\ln {17}}{16\ln 2}$ but less than
$\frac{\ln 3}{2\ln 2}$. \slshape Is it true for the further Fermat
primes the relation \upshape

$$
S_p(2^p x)=2p S_p(2^{\frac{p+1}{2}}x)-p S_p(2x) ?
$$

Unfortunately, our method (\cite{4}) for receiving such relations is
too tiring and does not give anything for large modulos. It is very
interesting to understand if the cases when 2 is (semi)primitive
root of $p$ are all the cases when we have a binomial relation of
the form $S_p(2^p x)= ap S_p(2x)?$

\end{document}